%% file: JiaZhaoTong14TSG_final.tex
\documentclass[journal]{IEEEtran}

\usepackage{graphicx,epsf,psfrag,amsmath,amssymb,amsfonts,verbatim}
\usepackage[mathscr]{eucal}
\usepackage{algorithmic}
\usepackage{verbatim,pifont,color,epsf,psfrag,amssymb,amsfonts,amsmath}
\usepackage{caption}
\usepackage{subcaption}

\IEEEoverridecommandlockouts                              

\input LTmacros

\def\DA{\mbox{\tiny DA}}
\def\RT{\mbox{\tiny RT}}

\begin{document}
%
\title{An online learning approach to dynamic pricing \\
for demand response}
%
%
%

\author{{\large Liyan Jia, Lang Tong and Qing Zhao}
\thanks{This work was supported by the National Science Foundation under
Grant CNS-1135844 and CNS-1248079 and by the Army Research Office under Grant W911NF-12-1-0271.
 }
\thanks{Liyan Jia and Lang Tong are  with School of Electrical and Computer Engineering, Cornell University, Ithaca, NY, USA, 14850. Email: {\tt \{lj92,lt35\}@cornell.edu}.}%
\thanks{Qing Zhao is with School of Electrical and Computer Engineering, University of California, Davis, CA, USA, 95616.  Email: {\tt qzhao@ucdavis.edu}.}%
\thanks{Part of this work was presented at the 51st Annual Allerton Conference on Communication, Control, and Computing, Oct. 2013.
}%
}

\markboth{Submitted to IEEE Transactions on Smart Grid}%
{Submitted to IEEE Transactions on Smart Grid}
%



\maketitle

\begin{abstract}
In this paper, the problem of optimal dynamic pricing for retail electricity with an unknown demand model is considered. Under the day-ahead dynamic pricing (a.k.a. real time pricing) mechanism, a retailer obtains electricity in a two-settlement wholesale market and serves its customers in real time. Without knowledge on the aggregated demand function of its customers, the retailer aims to maximize its retail surplus by sequentially adjusting its price based on the behavior of its customers in the past. An online learning algorithm, referred to as piecewise linear stochastic approximation (PWLSA), is proposed.  It is shown that PWLSA achieves the optimal rate of learning defined by the growth rate of cumulative regret.  In particular, the regret of PWLSA is shown to grow logarithmically with respect to the learning horizon, and no other on-line learning algorithm can have the growth rate slower than that of PWLSA.  Simulation studies are presented using traces of actual day-ahead prices, and PWLSA compares favorably under both static and dynamically changing parameters.

\end{abstract}


\section{Introduction}

As a key feature of a future smart grid, demand response is an effective way to improve power system operation efficiency, hedge the risk of energy supply shortage, and enhance social welfare. Based on characteristics of the interaction between a retailer and its consumers, demand response can be classified into two categories: demand responses with direct control and ones with indirect control; the former refers to programs in which the consumers enjoy lower electricity rate by allowing the retailer to shed load in case of emergency. The latter refers to approaches of influencing the consumers' consumption through dynamic pricing of electricity. In this paper, we focus on the latter.

We assume that the retailer employs a real-time pricing mechanism, referred to as day-ahead dynamic price (DADP), under which the retailer posts the hourly prices of electricity one day ahead. First proposed by Borenstein, Jaske, and Rosenfield \cite{Borenstein&etc:02} and referred to as the real-time pricing (RTP), DADP has been implemented in practice \cite{Borenstein&etc:02, Hopper&Goldman&Neenan:05EJ}. A key advantage of DADP is that a customer has the short-term price certainty with which it can optimize its consumption.

We assume that the retailer obtains electricity from a commonly adopted a two-settlement wholesale market, consisting of a day-ahead market and real-time market. In the day-ahead market, both the generators and retailers offer bids for the next day. Based on the submitted bids, the system operator schedules the day-ahead dispatch and clears the market with the day-ahead price. In real-time operations, the system operator adjusts the day-ahead dispatch according to the actual operation condition and sends dispatch signal to all participants to maintain the system balance. The amount of electricity deviated from the day-ahead schedule is settled according to the real-time price.


If the retailer knows how its customers respond to the retail price through their individual demand functions, it can choose the price to optimize a particular objective, \eg
the social welfare or its own profit subject to regulations. Obtaining the demand functions of its customers, however, is nontrivial because a customer is likely to consider such information private; neither the willingness of sharing nor the correctness of the shared information can be assumed.

In this work, we focus on optimal dynamic pricing under unknown demand functions. We take an online learning approach where the retailer learns the behavior of its customers by observing their response to carefully designed prices.  The basic principle of online learning is to achieve a tradeoff between  ``exploration'' and ``exploitation;'' the former represents the need of using sufficiently rich pricing signals to achieve learning accuracy, whereas the latter stands for the need of capturing as much reward as possible based on what has been learned.

In the classical online learning theory, the performance of a learning algorithm is measured by the notion of cumulative regret.  For the pricing problem at hand, the regret is defined as the difference between the retail surplus associated with the actual aggregated demand function and the surplus achieved by an online learning algorithm. While the cumulative regret $R_T$ grows with the learning horizon $T$, the rate of growth, $R_T/T$, of a well designed on-line learning algorithm typically diminishes, which implies that, for the infinite horizon problem, the profit achieved per unit time without knowing the demand function matches that when the demand function is known.  Therefore, a relevant performance metric is the growth rate of regret $R_T$ vs. $T$.

\subsection{Summary of results}

The basic problem setting involves two players: a retailer (an electricity distributor or aggregator) who offers its customer day-ahead hourly dynamic prices and its customers with price responsive demands.  We focus on the case when the customer demands are elastic and can be described by a random affine model, which arises naturally for thermal control applications.

The main result of this paper is twofold.  First,  under the DADP mechanism, we propose a simple online learning algorithm, referred to as piecewise linear stochastic approximation (PWLSA), that has the logarithmic rate of growth in regret, i.e. $R_T(T)=\Theta(\log T)$.

On the other hand, we show that no other on-line learning algorithm can have the rate slower than that of PWLSA. Thus PWLSA is order optimal. To achieve the optimal rate of learning, we deviate the standard on-line learning approach by first analyzing the mechanism of the two-settlement wholesale electricity market and calculate the retail surplus of the retailer as a wholesale market participant in a simple set-up. The result shows that the retailer's loss of surplus is proportional to the 2-norm deviation of the real-time consumption from the day-ahead schedule.


To demonstrate the learning performance, we also conduct simulations to compare PWLSA with the Greedy Method based on the actual data. In both cases with static and dynamically changing parameters of the demand model, PWLSA outperformed the greedy method and converged fast towards the optimal price.


\subsection{Related work}

The problem of dynamic pricing for demand response assuming known demand functions has been extensively studied.  See, for example, \cite{Borenstein&etc:02, Borenstein:05,Hopper&Goldman&Neenan:05EJ}, which adopted a similar pricing scheme as considered in this paper and \cite{Carrion&Etal:07TPS, Conejo&Etal:08TPS,Yang&Etal:12TPS} for more general settings. A precursor of the work presented here is \cite{Jia&Tong:12Allerton} where a parametric form of demand function was obtained.  In \cite{Jia&Tong:13CDC}, the tradeoff between retail profit and consumer surplus was characterized under a Stackelberg formulation with known demand functions.

The general problem of online learning for dynamic pricing has been studied extensively in multiple communities. This problem can be formulated as a multi-armed bandit (MAB) problem by treating each possible price as an arm. When the price can only take finite possible values, the problem becomes the classic MAB for which Lai and Robbins showed that the optimal regret growth rate is $\Theta(\log T)$ when the arms generate independent reward \cite{Lai&Robbins:85AAM}. When the price takes value from an uncountable set, the dynamic pricing problem is an example of the so-called continuum-armed bandit introduced by Agrawal in \cite{Agrawal:95SIAM} where the arms form a compact subset of $\Rc$. An online learning policy with regret order of $O(T^{3/4})$ was proposed in \cite{Agrawal:95SIAM} for any reward function satisfying Lipschitz continuity. Further development on the continuum-armed bandit under various assumptions of the unknown reward function can be found in \cite{Kleinberg:04,Auer&etal:07,Cope:09}.   The reason that PWLSA proposed in this paper achieves a much better regret order ($\Theta(\log T)$) than in the case of a general continuum-armed bandit is due to the specific linearly parameterized demand which leads to a specific quadratic cost/reward function. A similar message can be found in \cite{KleinbergLeighton:03,RusmevichientongTsitsiklis:10,BroderRusmevichientong:12,Harrison&etal:SC11,Zhai&etal:Asilomar11} where different regret orders were shown to be achievable under different classes of demand models for dynamic pricing.

The problem considered in this paper deals with linearly parameterized demand functions, thanks to the closed-form characterization of the optimized demand function for thermal dynamic load obtained in \cite{Jia&Tong:13CDC}. The learning approach proposed in this paper is rooted from a stochastic approximation problem originally formulated by Lai and Robbins \cite{Lai&Robbins:79AS,Lai&Robbins:82AAM} where the authors considered a form of optimal control problem when the model contains unknown parameters and the cost of control is explicitly modeled. For scaler models, Lai and Robbins showed in \cite{Lai&Robbins:79AS,Lai&Robbins:82AAM} that the cumulative regret (if translated from our definition) of a simple linear stochastic approximation scheme grows at the rate of $O(\log T)$. However, it is not clear whether such growth rate is the lowest possible.  Our result provides a generalization to the vector case with a lower bound for general policies. In addition, our approach also allows the consumers to have variable demand levels whereas the algorithm presented in  \cite{Lai&Robbins:79AS,Lai&Robbins:82AAM} only allows a single constant demand target.

Also related is the work of Bertsimas and Perakis \cite{Bersimas&Perakis:06} who tackled the problem as a dynamic program with incomplete state information. The authors showed in numerical simulations that considerable gain can be realized over the myopic policy where the price in the next stage is based on the least squares estimate of the model parameter.  When the parameters are assumed to be random, Lobo and Boyd considered the same problem under a Bayesian setting \cite{Lobo&Boyd:03Informs} and proposed a randomized policy via a dithering mechanism. In both cases, the rate of learning is not characterized.

Machine learning techniques have been applied to pricing problems in electricity markets, although there seems to be limited literature on discovering real-time price with unknown demand functions at the retail level.  While such problems can be viewed as part of the general learning problem discussed above, the nature of electricity market and electricity demand impose special constraints. When the market has multiple strategic generators, Garcia et al. proposed an online learning algorithm which converges to the Markov perfect equilibria \cite{Garcia&etal:05OP}. A related learning problem of bidding strategy of a retailer in the wholesale market when the supply functions of the generators are unknown has been studied.  See \cite{RahimiKianSadeghiThomas05PES,QiuPeetersDeconinck09ISAP,Pinto&Etal:11ISAP} where Q-learning techniques have been applied. Some other research focuses on developing learning methods for optimal demand response. See \cite{Taylor&Mathieu:13STPS} for index policy by formulating the demand control as a restless bandit problem, and \cite{O'Neill&etal:10SGC} for a reinforcement learning solution to a partially observable Markov decision process (MDP) problem.

\section{Structure of Wholesale Electricity Market }
\label{sec:whole}
In this section, we discuss a two-settlement system stylized from the deregulated wholesale market in the United States. The market consists of a day-ahead market and a real-time market. The day-ahead market serves as a planning mechanism for participants, and its settlement is financially binding. In the presence of uncertainties, the real-time market, on the other hand, addresses mismatches between the actual generation/consumption and that planned in the day-ahead market. In the following discussion, we only consider the presence of retailers and generators, without other financial participants, such as virtual bidders.

As a participant in the two-settlement market, a retailer faces uncertainties in the wholesale market and that from the real-time consumptions of its customers. If the quantity of consumption is large, the retailer is not a price taker. Instead, its bidding curve and real-time purchase will affect the wholesale price. Using a simplified model, we argue in this section that it is to the retailer's benefit to match the real time consumption with the day-ahead dispatched value.  In particular, we motivate, by algebraic and economic arguments, that minimizing the 2-norm deviation of the real-time consumption maximizes the retail surplus.  This result motivates the specific form of the cost used in the regret definition in our online learning formulation of the problem.



\subsection{The day-ahead wholesale market}

In the day-ahead market, the independent system operator (ISO) schedules energy dispatch for the next day. Each electricity generator submits a cost curve $c(p)$ that represents the cost of serving $p$ units of electricity, while each retailer (or Load Serving Entity (LSE)) submits a utility curve $u(d)$ that models the benefit of getting served with $d$ units of electricity. Usually, the day-ahead market dispatch is calculated at the hourly time scale. Therefore, both the demand schedule $d$ and the generation schedule $p$ are $24$ dimensional vectors.

With all submitted offers and bids, the ISO solves an optimal power flow (OPF) problem to obtain the optimal dispatch under the objective of maximizing the social welfare. In its simplest form without complications of capacity constrained transmission networks and multiple participating agents, the OPF problem is of the following form,
\begin{equation}
  \label{eq:opf}
  \begin{array}{r l}
  \mbox{max}_{d,p} & u(d) - c(p) \\
  \mbox{s.t.} & d = p \\
  \end{array}
\end{equation}

The solutions, $d^{\DA}$ and $p^{\DA}$, represent the desired day-ahead dispatch of demand and generation. The day-head price is defined as the cost of serving next unit of energy. Therefore, it is the marginal cost of generating $p^{\DA}$, $\ie$ $\lambda^{\DA} = \frac{\partial c}{\partial p}(p^{\DA})$.

The clearing of the day-ahead market is financially binding in the sense that, regardless of the actual consumption in real time, the day-ahead payment from retailer to the system operator is settled as $(\lambda^{\DA})^{\mbox{\tiny T}} d^{\DA}$. The payment from the system operator to the generator is $(\lambda^{\DA})^{\mbox{\tiny T}} p^{\DA}$. Since the retailer's utility of using $d^{\DA}$ is $u(d^{\DA})$, the retail surplus is
\begin{equation}
S^{\DA}_{\mbox{\tiny retail}} = u(d^{\DA}) - (\lambda^{\DA})^{\mbox{\tiny T}} d^{\DA}.
\end{equation}

\subsection{The real-time wholesale market}
\label{ssec:RTmarket}
The actual consumption and generation in real time $d^{\RT}$ and $p^{\RT}$, however, are nominally different from the day ahead dispatch. Consequently, the real-time price will deviate from the day-ahead price.  In particular,  if the cost function of generation in real-time is $\tilde{c}(p)$, the real-time price is calculated as $\lambda^{\RT} = \frac{\partial \tilde{c}}{\partial p}(d^{\RT})$, which stands for the cost of serving the next unit of electricity in real time.

Different from the day-ahead settlement, the real-time settlement only applies to the difference between the day-ahead schedule and the real-time consumption. This means that the payment from the retailer to the system operator is $(\lambda^{\RT})^{\mbox{\tiny T}} (d^{\RT} - d^{\DA})$ if positive. Otherwise, this quantity represents the compensation from the system operator to the retailer.


Therefore, if the real-time consumption matches the day-ahead dispatch, there is no real-time payment. The total retail surplus is still $S^{\DA}_{\mbox{\tiny retail}}$. If the actual consumption $d^{\RT}$ is different from $d^{\DA}$, the retail surplus is
\begin{equation}
S^{\RT}_{\mbox{\tiny retail}} = u(d^{\RT}) - [(\lambda^{\DA})^{\mbox{\tiny T}} d^{\DA} + (\lambda^{\RT})^{\mbox{\tiny T}} (d^{\RT} - d^{\DA})],
\end{equation}
where the first term is the utility of the retailer from delivering $d^{RT}$ to its consumer, and the second term is the total payment to the wholesale market. Therefore, the surplus loss due to deviation of $d^{\RT}$ from $d^{\DA}$ is
\begin{equation}
\Delta S_{\mbox{\tiny retail}} = S^{\DA}_{\mbox{\tiny retail}} - S^{\RT}_{\mbox{\tiny retail}}.
\end{equation}




Based on the Taylor expansion of $u(d^{\RT})$, we can approximate $\Delta S_{\mbox{\tiny retail}} $ as shown in the following lemma. The complete proof is included in the Appendix.

\begin{lemma}
\label{lemma:square}
Under the assumption that cost c(p) and utility u(d) are twice differentiable,
\begin{equation}
\Delta S_{\mbox{\tiny retail}} \approx \theta(d^{\RT} - d^{\DA})^{\mbox{\tiny T}} (d^{\RT} - d^{\DA}),
\end{equation}
where $\theta$ is a constant independent of $d^{\RT}$ and $d^{\DA}$
\end{lemma}

Therefore, the objective of maximizing retail surplus is equivalent to minimizing the squared deviation of the real-time demand to the day-ahead dispatch.

\begin{center}
\begin{figure}[!ht]
\begin{center}
\begin{psfrags}
\psfrag{P0}[c]{\small{$\lambda^{\DA}$}}
\psfrag{P1}[c]{\small{$\lambda^{\RT}$}}
\psfrag{D0}[c]{\small{$d^{\DA}$}}
\psfrag{D1}[c]{\small{$d^{\RT}$}}
\psfrag{O}[c]{}
\psfrag{Price}[l]{\small{price}}
\psfrag{Quantity}[l]{\small{Quantity}}
\psfrag{Area I}[l]{\small{Area I}}
\psfrag{Area II}[l]{\small{Area II}}
\psfrag{Area III}[l]{\small{Area III}}
\psfrag{RT Supply}[l]{\small{real-time supply function}}
\psfrag{DA Supply}[l]{\small{day-ahead supply function}}
\psfrag{Demand function}[l]{\small{demand function}}
\includegraphics[width = 2.0in]{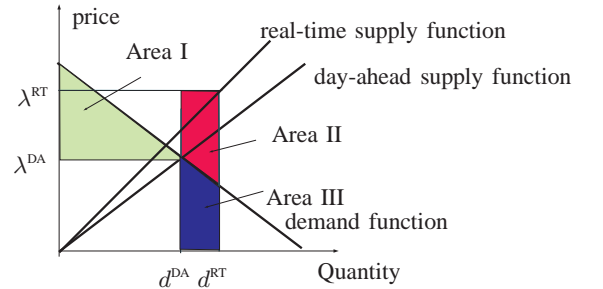}~~~%
\end{psfrags}
\end{center}
\caption{Real-time market equilibrium}
\label{fig:RTeq}
\end{figure}
\end{center}

The result above can also be illustrated in the Price-Quantity plane as shown in Fig.~\ref{fig:RTeq}. The demand function presents the optimal quantity of energy required from the retailer given the price. It is actually the derivative of the utility function $u(d)$. The area below the line is the integration, which is exactly the utility value with quantity $d$. Similarly, the day-ahead and real-time supply function stand for the optimal quantity of generation to the generator if the price is given. The crossing point $(d^{DA}, \lambda^{\DA})$ is the day-ahead equilibrium, the same as calculated from (\ref{eq:opf}).  Subtracting the day-ahead payment from the utility,  Area I represents the day-ahead retail surplus.

In the real-time market, the real-time consumption $d^{\RT}$ deviates from $d^{\DA}$, and the real-time price, $\lambda^{\RT}$, is determined by the real-time supply function. Area III is the additional utility gained by consuming $d^{\RT}$, while the sum of area II and III is the real time payment. Therefore, Area II represents the retail surplus loss, and the loss grows in the order of $||d^{\RT}-d^{\DA}||_2^2$---the 2-norm deviation between the day-ahead scheduled consumption and the actual real-time consumption.

\section{Price responsive demand}
\label{sec:demand}


In this section, we characterize the behavior of the consumers in a demand response program offered by a retailer that uses dynamic price to influence consumptions of its customers.  To this end, the retailer is particularly interested in parametric models of the demand function that captures the relation between price of electricity and the level of consumption.

Instead of assuming a particular model, we consider an  important engineering application for HVAC units\footnote{Heating, ventilation, and air conditioning units} based temperature control,  since it makes up most price responsive demand  \cite{EnergySurvey}. We establish in this section that the optimal demand response under the day ahead dynamic pricing is has an affine parametric form.

\subsection{Day-ahead dynamic pricing}

Following the conclusion in Section~\ref{sec:whole}, the retailer's optimal strategy is to minimize the difference between the real-time demand and the day-ahead dispatch. In this paper, we assume that the retailer is to influence the aggregated consumption of its customers via retail pricing. Specifically, we consider a specific form of dynamic pricing, Day-Ahead Dynamic Pricing (DADP), also referred to as Real-Time Pricing (RTP) by Borenstein et al. \cite{Borenstein&etc:02}.


DADP works in the following way. In the day-ahead market, the retailer offers the consumer hourly retail price $\pi$ for the day ahead, according to the prediction of the wholesale market and projected demand response. In real time, knowing the entire price trajectory for the day, a consumer dynamically determines her own energy consumption. The payment from a consumer to the retailer is settled as the product of day-ahead price and real-time consumption. The retailer meets aggregated demand by purchasing electricity from the wholesale market.

\subsection{Optimal demand respose}

Given DADP posted one day ahead, the consumer optimizes her energy consumption.  Here we assume a general linear dynamic model that captures the relation between the consumer utility and her consumption of electricity.  This general model arises from the classical model for HVAC based temperature control where thermal storage is involved \cite{EnergySurvey,Bargiotas&Birddwell:88ITPD} and is shown to be reasonably accurate for residential temperature control \cite{Yu&etal:13TSG}.



For consumer $k$, let $x^{(k)} = (x^{(k)}_1,...,x^{(k)}_{24})$ and $a = (a_1,...,a_{24})$ denote the average indoor temperature and outdoor temperature in each hour, respectively. The Birdwell model of HVAC is given by
\begin{equation}
x^{(k)}_i=x^{(k)}_{i-1}+ \alpha^{(k)} (a_i-x^{(k)}_{i-1}) - \beta^{(k)} u^{k}_{i} + \xi^{(k)}_i,
\label{eq:hvacmodel}
\end{equation}
where  $u^{(k)} = (u^{(k)}_1,...,u^{(k)}_{24})$ is the vector of control variable representing the total amount of electricity drawn by the HVAC unit during each hour and $\xi^{(k)} = (\xi^{(k)}_1,...,\xi^{(k)}_{24})$ the process noise. System parameters $\alpha^{(k)}$ ($0<\alpha<1$) and $\beta^{(k)}$ model the insolation of the building and the efficiency of the HVAC unit in consumer $k$'s house. Note that the above equation applies to both heating and cooling scenarios but not simultaneously. We focus herein the cooling scenario ($\beta^{(k)} > 0$) and the results apply to heating ($\beta^{(k)} < 0$) as well.

Using a linear combination of total cost and squared deviation of indoor temperature from desired temperature as the objective function to minimize, the consumer $k$'s energy consumption can be modeled by the following stochastic optimization problem,
\begin{equation}
\begin{array}{r l}
\min & \mathbb{E}\left\{\sum_{i=1}^{24} [ \kappa (x^{(k)}_i - \xi^{(k)}_{i})^2] + \pi^{\tiny{\text{T}}} u^{(k)} \right\} \\
\mbox{s.t.} & x^{(k)}_{i} = x^{(k)}_{i-1} + \alpha^{(k)} (a_{i} - x^{(k)}_{i-1}) - \beta^{(k)} u^{(k)}_i + \xi^{(k)}_i, \\
& y^{(k)}_i=(x^{(k)}_i,a_i)+ \nu^{(k)}_i, \\
\end{array}
\end{equation}
where $y^{(k)}=(y^{(k)}_1,...,y^{(k)}_{24})$ is the observation vector, $\nu^{(k)}=(\nu^{(k)}_1,...,\nu^{(k)}_{24})$ the observation noise vector, $\kappa$ the weight factor and $\xi_i$ the desired temperature for hour $i$.

The solution of the above stochastic optimization can be obtained in closed form  via direct backward induction.  More significantly, it is shown in
\cite{Jia&Tong:13CDC}  that, after aggregation over all consumers, the total demand is an affine function of the retail price.

\begin{theorem}[\cite{Jia&Tong:13CDC}]
\label{thm:opt_demand}
Assume that the process noise $\xi^{(k)}$ and $\nu^{(k)}$ are Gaussian distributed with zero mean for each consumer $k$. With the fixed retail price $\pi$, the optimal aggregated residential demand response has the following matrix form and properties,
\begin{equation}
d^{\RT}= \sum_k u^{(k)} =b -A\pi+w,
\label{eq:demand}
\end{equation}
where  the factor matrix $A$ is positive definite, $b$ and $A$ are deterministic, depending only on the dynamic system parameters, and $w$ is a random vector with zero mean.
\end{theorem}

\section{Dynamic Retail Pricing via Online Learning}
\label{sec:pricing}
\subsection{Pricing policy and regret}

As discussed in Section~\ref{sec:whole}, minimizing the demand side surplus loss is equivalent, approximately, to minimizing the squared deviation of real-time electricity consumption, $d^{\RT}$, from the day-ahead optimal dispatch, $d^{\DA}$.

Formally, define the $t$-th day's expected surplus loss as the 2-norm of the deviation of the real-time consumption from the day-ahead dispatch, $\ie$ $L_t \defeq \mathbb{E}[||d_t^{\RT} - d^{\DA}_t||_2^2]$, where $d^{\DA}_t$ and $d^{\RT}_t$ are the day-ahead and real-time demands for day $t$.

Assuming the linear demand function in Theorem~\ref{thm:opt_demand}, for the purpose of obtaining a performance upper bound, we consider the case that the parameters in (\ref{eq:demand}), $A$ and $b$, are known to the retailer. At day $t$, the optimal retail price is given by
\begin{equation}
\pi^*_t = \mbox{arg} \min_{\pi_t} \mathbb{E}[||d^{\RT}_t - d^{\DA}_t||_2^2] = A^{-1}(b - d^{\DA}_t),
\label{eq:optpi}
\end{equation}
and the corresponding minimum surplus loss is only caused by the exogenous random fluctuations (such as the outdoor temperature).  Specifically, the minimized expected loss is
\begin{equation}
\mathbb{E}[||b - A\pi^*_t + w_t - d^{\DA}_t||_2^2] = ||\Sigma_w||_2,
\end{equation}
where $\Sigma_w$ is the covariance matrix of demand model noise $w$ in (\ref{eq:demand}). Notice that the minimized surplus loss is independent of the day-ahead dispatch $d^{\DA}_t$.

However, it is nontrivial for the retailer to obtain the exact parameters of the demand functions of its customers because a customer is likely to consider such information private. At day $t$, the only information available to the retailer is the record of previous electricity consumption up to $t-1$ and day-ahead dispatch up to $t$.  Formally, the retail pricing policy is defined as follows,

\begin{definition}
The retail pricing policy $\mu =(\mu_t)$ is a sequence of mappings where
$\mu_t$ maps the consumption history and day-ahead demand dispatch to the price vector of day $t$. In particular,  letting $\pi_t^{\mu}$ be the price vector under policy $\mu$, we have
 \begin{equation}
  \pi_{t}^{\mu} = \mu_t (d_{0}^{\RT},...,d_{t-1}^{\RT}, d_0^{\DA},...,d^{\DA}_{t-1}, d^{\DA}_{t}),
  \label{eq:policy}
  \end{equation}
  where $d_{i}^{\DA},$ and $d^{\RT}_i$ are the day-ahead dispatch and real-time electricity consumption for day $i$.
  \hfill $\square$
\end{definition}

As for a particular policy $\mu$, the regret $R_t^{\mu}$ at day $t$ is defined as the increase of surplus loss compared with using the optimal price, $\pi^*_t$, which means that
  \begin{equation}
  \label{eq:regret}
  \begin{array}{r l}
  R_t^{\mu} & \defeq \mathbb{E}[||b - A\pi_t^{\mu} + w_t - d^{\DA}_{t}||^2_2 - ||\Sigma_w||_2] \\[0.2em]
  & = \mathbb{E}[||b - A\pi_t^{\mu} - d_t^{\DA}||^2_2].
  \end{array}
  \end{equation}
Because maximizing the surplus is equivalent to minimizing the regret, we'll focus next on the increasing rate of the cumulative regret up to day $T$, $\sum_{t=1}^T R_t^{\mu}$.

\subsection{Lower bound on the growth rate of regret}

To gain insights into the type of lower bound on the regret, we consider first a simple example when part of the parameters are known.  Intuitively, the advantage of knowing partially the parameters should lead to a lower growth rate of regret.

In particular, recall the stochastic affine demand function (\ref{eq:demand}) where we assume parameter $A$ is known but $b$ is unknown. Consider the following dynamic pricing policy, $\tilde{\mu}$, given by
\begin{equation}
\label{eq:knownA}
\pi^{\tilde{\mu}}_{t} = \bar{\pi}_{t-1} + A^{-1}(\bar{d}^{\RT}_{t-1} - d_{t}^{\DA}),
\end{equation}
where $\bar{\pi}_{t-1}$ and $\bar{d}_{t-1}$ are the average price and demand up to day $t-1$, $\ie$
\begin{equation}
\bar{\pi}_{t-1} = \frac{1}{t}\sum_{i=0}^{t-1}\pi_i \text{ and, } \bar{d}^{\RT}_{t-1} = \frac{1}{t}\sum_{i=0}^{t-1}d_i^{\RT}.
\end{equation}
According to (\ref{eq:regret}), straight forward calculation gives that the regret for day $t$ is
\begin{equation}
R_t^{\tilde{\mu}} = \mathbb{E} [||\frac{1}{t}\sum_{i=0}^{t-1} w_i ||_2^2] = \frac{1}{t} ||\Sigma_{w}||_2.
\end{equation}
Therefore, the aggregated regret,
\begin{equation}
\sum_{t=1}^{T} R_t^{\tilde{\mu}} = \sum_{t=1}^{T} \frac{1}{t} ||\Sigma_{w}||_2 \le (1+\log T)||\Sigma_{w}||_2
\end{equation}

Therefore, with the knowledge of the demand function parameter $A$, the policy $\tilde{\mu}$ achieves the aggregated regret $O(\log T)$ for any $b$ and any arbitrary sequence of $\{d_t^{\DA}\}$.

To establish the actual lower bound on the growth rate of regret, we formulate a game that, after the retailer proposes a deterministic pricing policy $\mu$, there exists an adversary designing parameters of the demand function. The adversary is to create the worst loss to the retailer while the retailer tries to minimize the largest possible loss. In other words, we consider the following the min-max regret as the objective,
\[
    \min_{\mu} \max_{b,A} \sum_{t=1}^{T} R_t^{\mu}.
\]

The following theorem shows that in the min-max sense, the growing rate of the cumulative regret can not be lower than $\log T$.
\begin{theorem}
For any pricing policy $\mu$ as defined in (\ref{eq:policy}), there exist some $(A, b, d_0^{\DA},...,d^{\DA}_{t-1}, d^{\DA}_{t},...)$ to make the cumulative regret, $\sum_{t=1}^{T} R_t^{\mu}$, grows at least at the rate of $\log T$.
  \label{thm:lowerbound}
\end{theorem}

Proof: see the Appendix.


\subsection{PWLSA: a rate optimal learning policy}


In this section, we propose a policy that achieves the lower bound on the regret growth rate; it is thus optimal in the sense of having the lowest rate of growth.  Referred to as piecewise linear stochastic approximation (PWLSA) policy, the proposed policy is an extension of the stochastic approximation approach of Lai and Robin \cite{Lai&Robbins:82AAM} for scaler processes with a single desired optimal price.

If the day-ahead demand is the same for all days, stochastic approximation will use the previous average price as the nominal value and previous average demand as the feedback signal to calculate the next price, as shown below,
\begin{equation}
\label{eq:singlelevel}
\pi^{\mbox{\tiny SA}}_{t} = \bar{\pi}_{t-1} + \gamma(\bar{d}^{\RT}_{t-1} - d^{\DA}),
\end{equation}
where $d^{\DA}$ is the constant day-ahead dispatch level, and the feedback factor $\gamma$ is a positive scalar.

For multiple day-ahead dispatch levels, we build adaptively a dictionary of day-ahead dispatch levels that have appeared before. Denote the dictionary at day $t$ as $\mathscr{D}_t$. For day $t+1$, if $d^{\DA}_{t+1} \in \mathscr{D}_t$, let $\mathscr{D}_{t+1} = \mathscr{D}_{t}$. Otherwise,  $\Dmsc_{t+1}=\Dmsc_t \bigcup \{d_{t+1}^{\DA}\}$. For each day-ahead dispatch level in $\mathscr{D} = \bigcup_{t=1}^{\infty} \mathscr{D}_t$, we keep a separate stochastic approximation to calculate the retail price, in a feedback control fashion similar as (\ref{eq:singlelevel}).

Therefore, for different $d_t^{\DA}$, we have a different linear function to calculate the next retail price. The policy is piecewise linear. Formally, the PWLSA policy,  $\mu^{\mbox{\tiny PWLSA}}$, is defined as,
\begin{definition}[PWLSA]
Assume for all $t \in \mathbb{N}^+$, $d_t^{\DA} \in \mathscr{D}$ and $\mathscr{D}$ is countable.
\begin{itemize}
  \item If $d_t^{\DA} \in \Dmsc_t$, then $\Dmsc_{t+1}=\Dmsc_t$ and
\begin{equation}
\pi_{t}^{\mbox{\tiny PWLSA}} =  \frac{1}{{|\mathscr{C}^{d_t^{\DA}}_t|}}\left(\sum_{k \in \mathscr{C}_{d_t^{\DA}}}  \pi_k^{\mbox{\tiny PWLSA}} + \gamma(d_k^{\mbox{\tiny PWLSA}} - d_{t}^{\DA})\right),
\end{equation}
where  $\mathscr{C}^{d_t^{\DA}}_t = \{k \in \mathbb{N}^+: k\le t-1, d_k^{\DA} = d_t^{\DA}\}$ and $|\mathscr{C}^{d_t^{\DA}}_t|$ is the total number of elements in $\mathscr{C}^{d_t^{\DA}}_t$.
 \item Otherwise, $d_t^{\DA} \notin \Dmsc_t$, then $\Dmsc_{t+1}=\Dmsc_t \bigcup \{d_t^{\DA}\}$ and
 \begin{equation}
 \pi_{t}^{\mbox{\tiny PWLSA}} = \tilde{\pi}_j,
 \end{equation}
 where $\tilde{\pi}_j$ is an arbitrary predetermined price.
\end{itemize}
\hfill $\square$
\end{definition}

The following theorem shows that PWLSA can achieve the optimal logarithmic regret order.
\begin{theorem}
\label{thm:logn}
Assume that day-ahead dispatch $d^{\DA}_t$'s are from a finite set, $\ie$ $|\mathscr{D}| < \infty$. If $\gamma \ge \frac{1}{2\lambda_{\min}(A)}$, where $\lambda_{\min}(A)$ is the minimum eigenvalue of $A$, then we have,
\begin{equation}
\sum_{t=1}^{T} R^{\mu^{\mbox{\tiny PWLSA}}}_t \sim O(\log (T)),
\end{equation}
\end{theorem}

Proof: see the Appendix.

Since $\log T$ is shown to be the optimal rate achievable, PWLSA is already the best in the sense of asymptotic growing rate of the regret. The conditions in the theorem are quite general. In practice, the cost functions from the generator and the utility functions from the retailer won't change often and are usually chosen from a few alternatives. Therefore, we can assume the total number of possible day-ahead dispatch levels, $|\mathscr{D}|$, to be finite. On the other hand, the consumers' demand function is from real data, which can be constrained by a compact set. Therefore, the bound of the minimum eigenvalue of $A$ is not hard to get with reasonable assumption.

\section{Simulation}
\label{sec:sim}

\subsection{Simulation set-up}
In this section, we conducted simulations based on the actual temperature records in Hartford, CT, from July 1st, 2012 to July 30th, 2012. The day-head price was also for the same period from ISO New England. The HVAC parameters for the simulation were set as: $\alpha = 0.5$, $\beta = 1$, $\mu=10$. The desired indoor temperature was set to be $18^{\circ}C$ for all hours. The size of aggregation was assumed to be 100.

\subsection{Learning static parameters}

First, we examined PWLSA's ability to identify the correct price if the parameters of the demand model remain the same. To make the comparison, we used the Greedy Method \cite{Anderson&Taylor:76,Lobo&Boyd:03Informs} as a benchmark. At each day, the Greedy Method makes maximum likelihood estimate of the parameters and uses the result as the correct parameters to calculate the ``optimal'' price.

\begin{figure}

        \centering
        \begin{subfigure}[b]{0.20\textwidth}
                \begin{psfrags}
                \psfrag{N}[c]{\tiny{Number of Days}}
                \psfrag{R}[c]{\tiny{Cumulative Regret}}
                \includegraphics[width=1.5in]{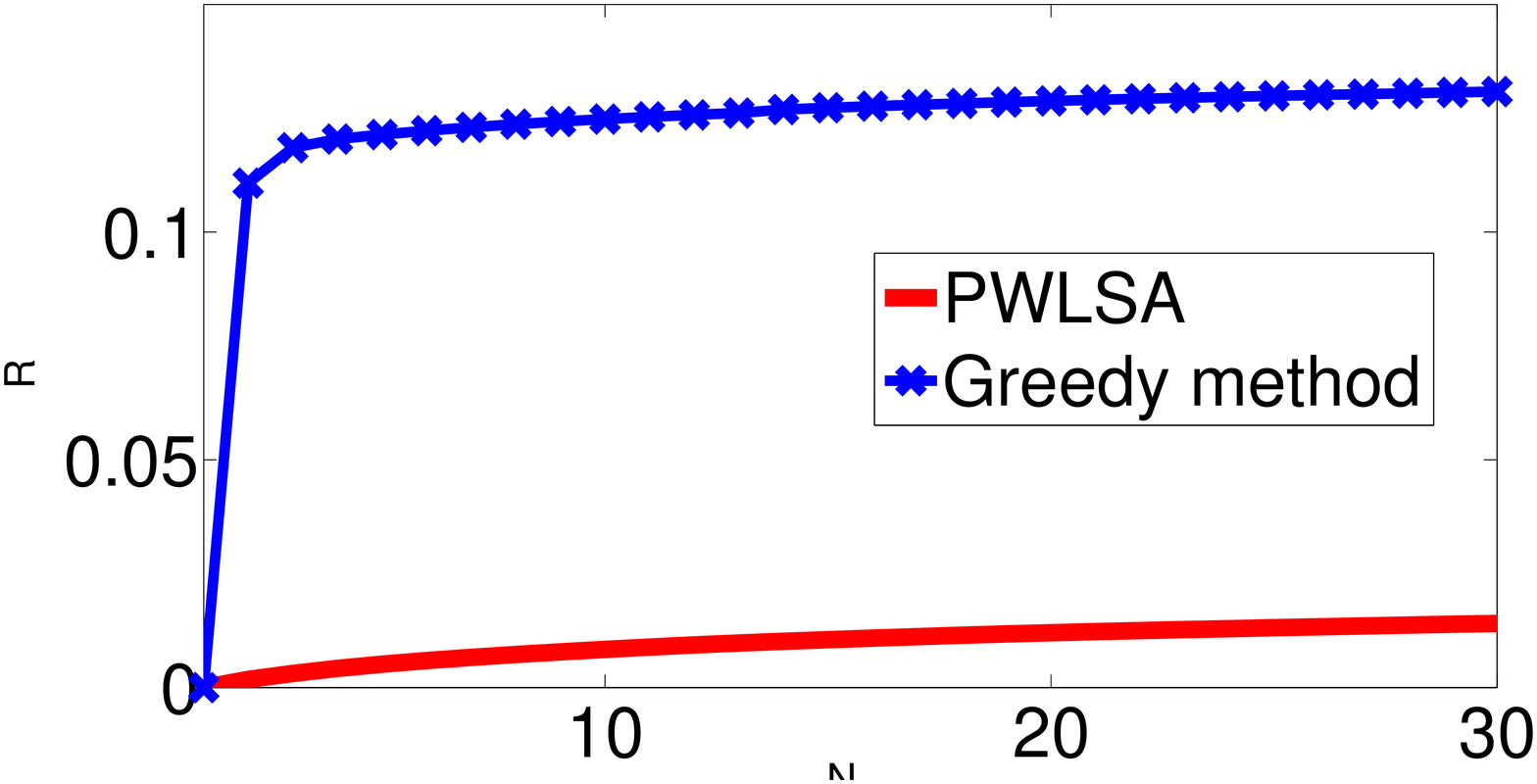}
                \caption{Cumulative regret}
                \label{fig:static_avg_regret}
                \end{psfrags}
        \end{subfigure}%
        \begin{subfigure}[b]{0.20\textwidth}
                \begin{psfrags}
                \psfrag{N}[c]{\tiny{Number of Days}}
                \psfrag{D}[c]{\tiny{Price Deviation}}
                \includegraphics[width=1.5in]{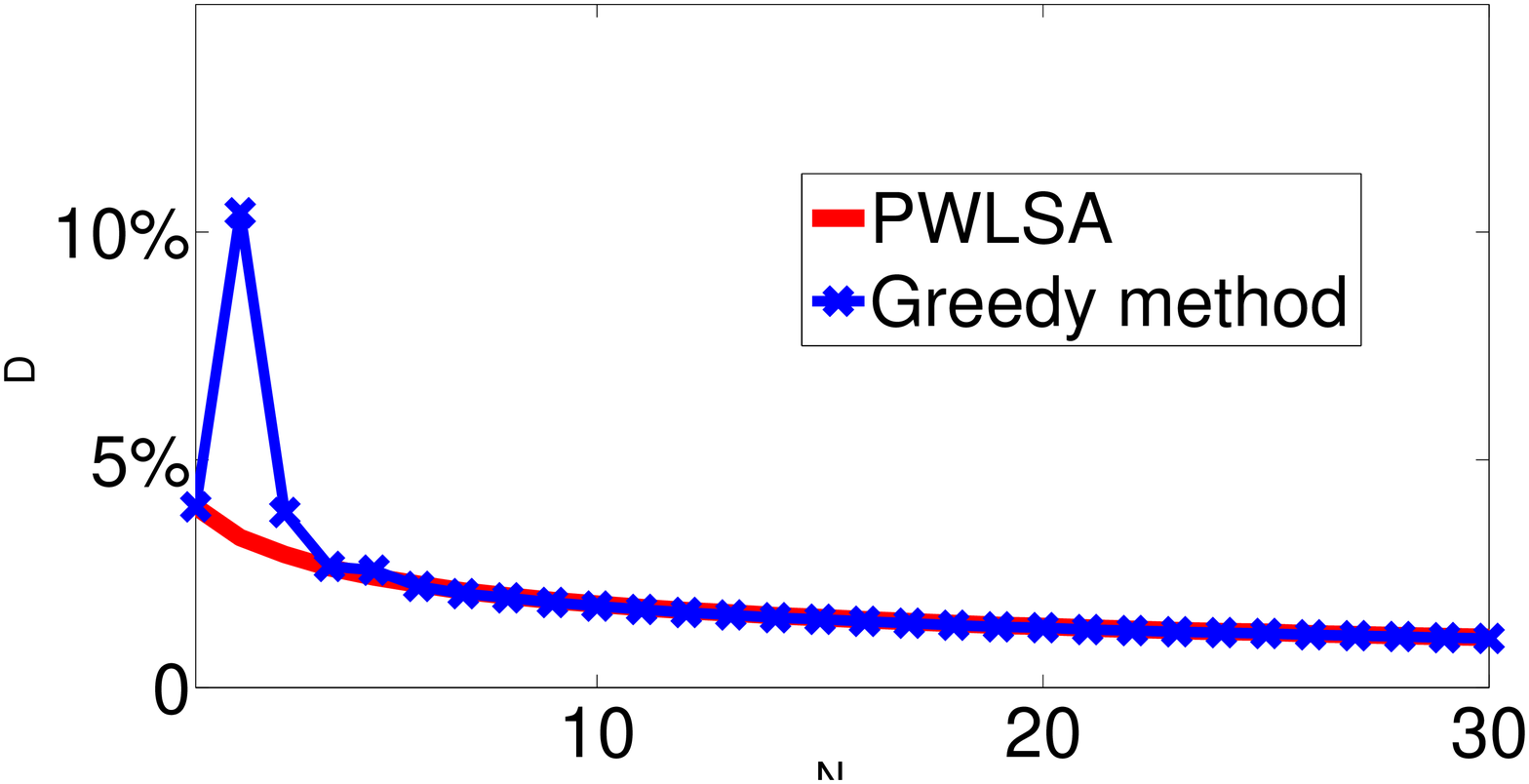}
                \caption{Price convergence}
                \label{fig:static_avg_price}
                \end{psfrags}
        \end{subfigure}%
\caption{Average performance comparison of PWLSA and the Greedy Method}
\label{fig:static_avg}
\end{figure}

Fig.~\ref{fig:static_avg} shows the average performance of PWLSA and Greedy Method over 10,000 Monte Carlo runs. In Fig.~\ref{fig:static_avg_regret}, the induced cumulative regrets of the two policies are compared. We could identify the logarithmic growth of the cumulative regret under PWLSA and significant cumulative regret increase by the Greedy Method. Fig.~\ref{fig:static_avg_price} shows the absolute percentage deviation of the prices under the two polices from the optimal price. We can see that Greedy Method performed extremely bad at the very beginning due to insufficient learning. After some days, the two policies both produced prices pretty close to the optimal one.

\begin{figure}

        \centering
        \begin{subfigure}[b]{0.21\textwidth}
                \begin{psfrags}
                \psfrag{N}[c]{\tiny{Number of Days}}
                \psfrag{R}[c]{\tiny{Price Ratio}}
                \includegraphics[width=1.5in]{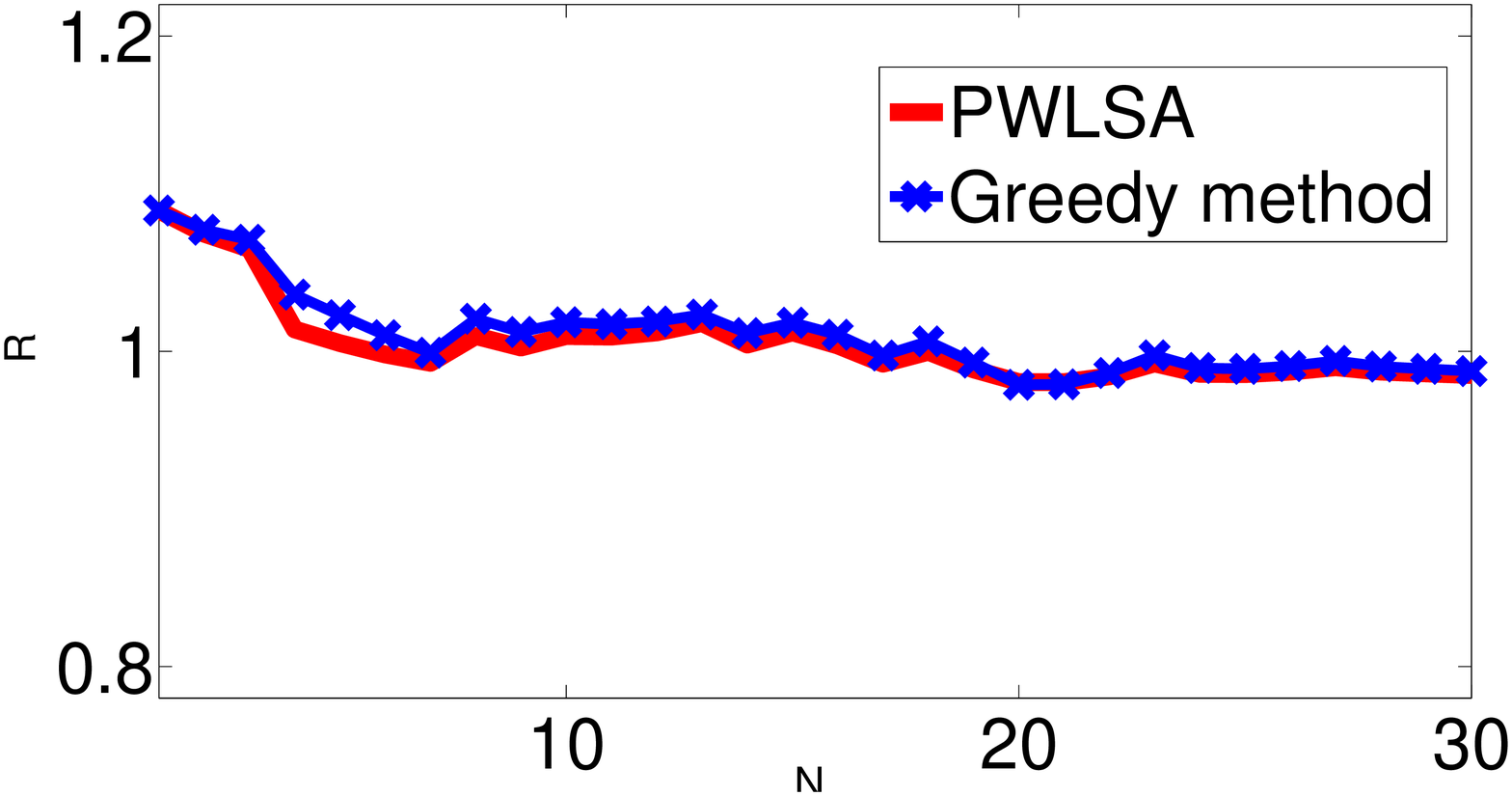}
                \caption{Scenario 1: The two had similar performance}
                \label{fig:static_scenario1}
                \end{psfrags}
        \end{subfigure}%
        \begin{subfigure}[b]{0.03\textwidth}
        \end{subfigure}
        \begin{subfigure}[b]{0.21\textwidth}
                \begin{psfrags}
                \psfrag{N}[c]{\tiny{Number of Days}}
                \psfrag{R}[c]{\tiny{Price Ratio}}
                \includegraphics[width=1.5in]{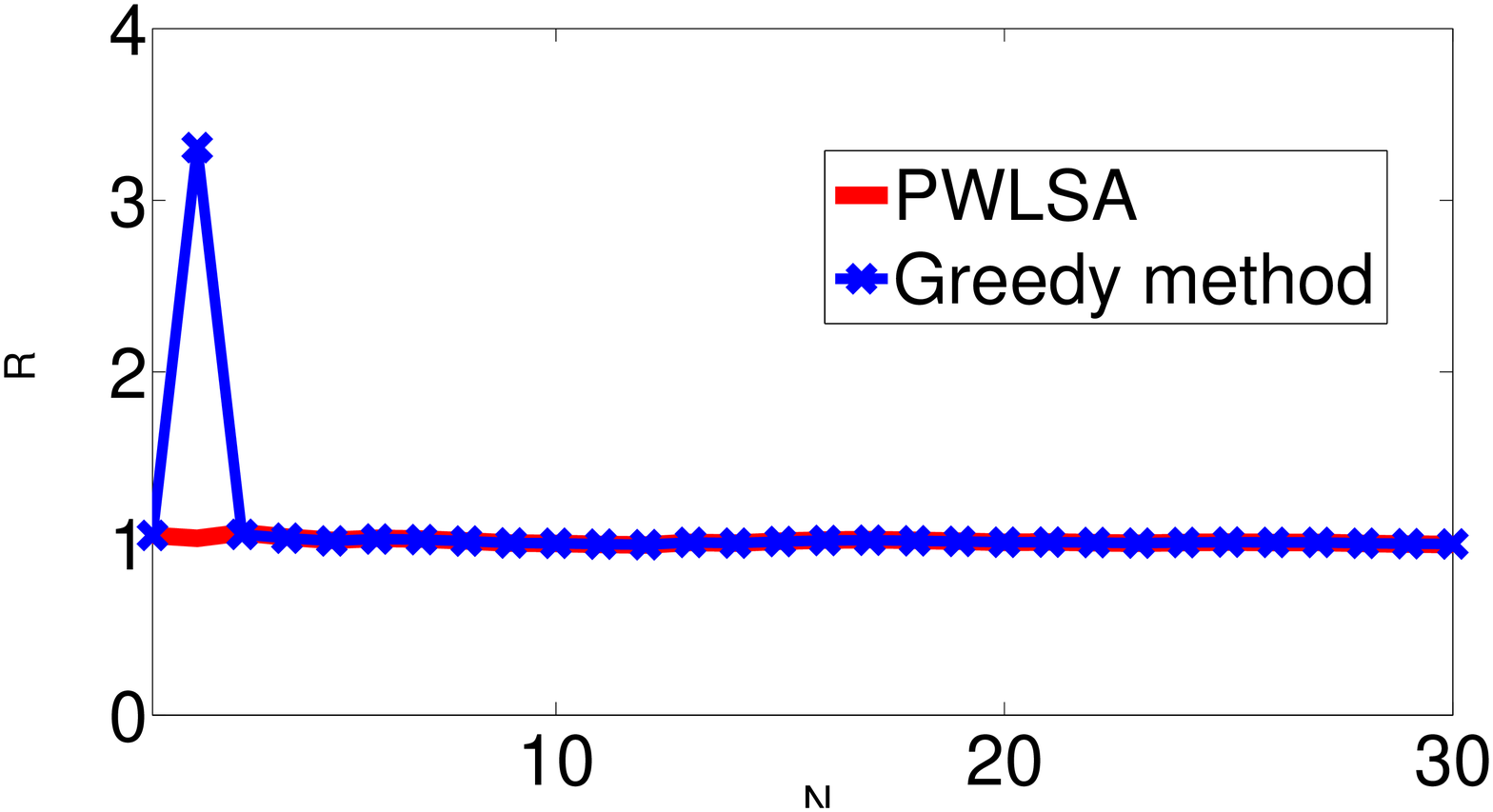}
                \caption{Scenario 2: the Greedy Method failed}
                \label{fig:static_scenario2}
                \end{psfrags}
        \end{subfigure}%
\caption{Scenario analysis of PWLSA and Greedy Method under static parameters}
\label{fig:static_scenario}
\end{figure}

After carefully investigating the simulated data, we found two typical scenarios as shown in Fig.~\ref{fig:static_scenario}. We used the ratio of the calculated price to the optimal price as y-axis, to show the fluctuation. In most of the cases as in Fig.~\ref{fig:static_scenario1}, the two polices gave similar performance and both converged to the optimal price fast. On the other hand, Fig.~\ref{fig:static_scenario2} shows one extreme scenario that Greedy Method run into the condition that is close to singularity, which leads to an abnormal price. Although this kind of scenarios happened rarely, it caused the wide performance gap between the Greedy Method and PWLSA.

\subsection{Learning dynamic parameters}
In the real world, the parameters of the demand model usually do not stay constant. They may follow some cycles or drifts. In this subsection, we tested the learning ability and robustness of PWLSA under dynamic unknown parameters. Besides the set of parameters above, we used $1.5A$ instead of $A$ to make the alternative set of parameters. We assumed the parameters followed a Markov Chain with these two sets as states. The transition probability to the other set was assumed to be $0.25$.

\begin{figure}
        \centering
        \begin{subfigure}[b]{0.20\textwidth}
                \begin{psfrags}
                \psfrag{N}[c]{\tiny{Number of Days}}
                \psfrag{R}[c]{\tiny{Cumulative Regret}}
                \includegraphics[width=1.5in]{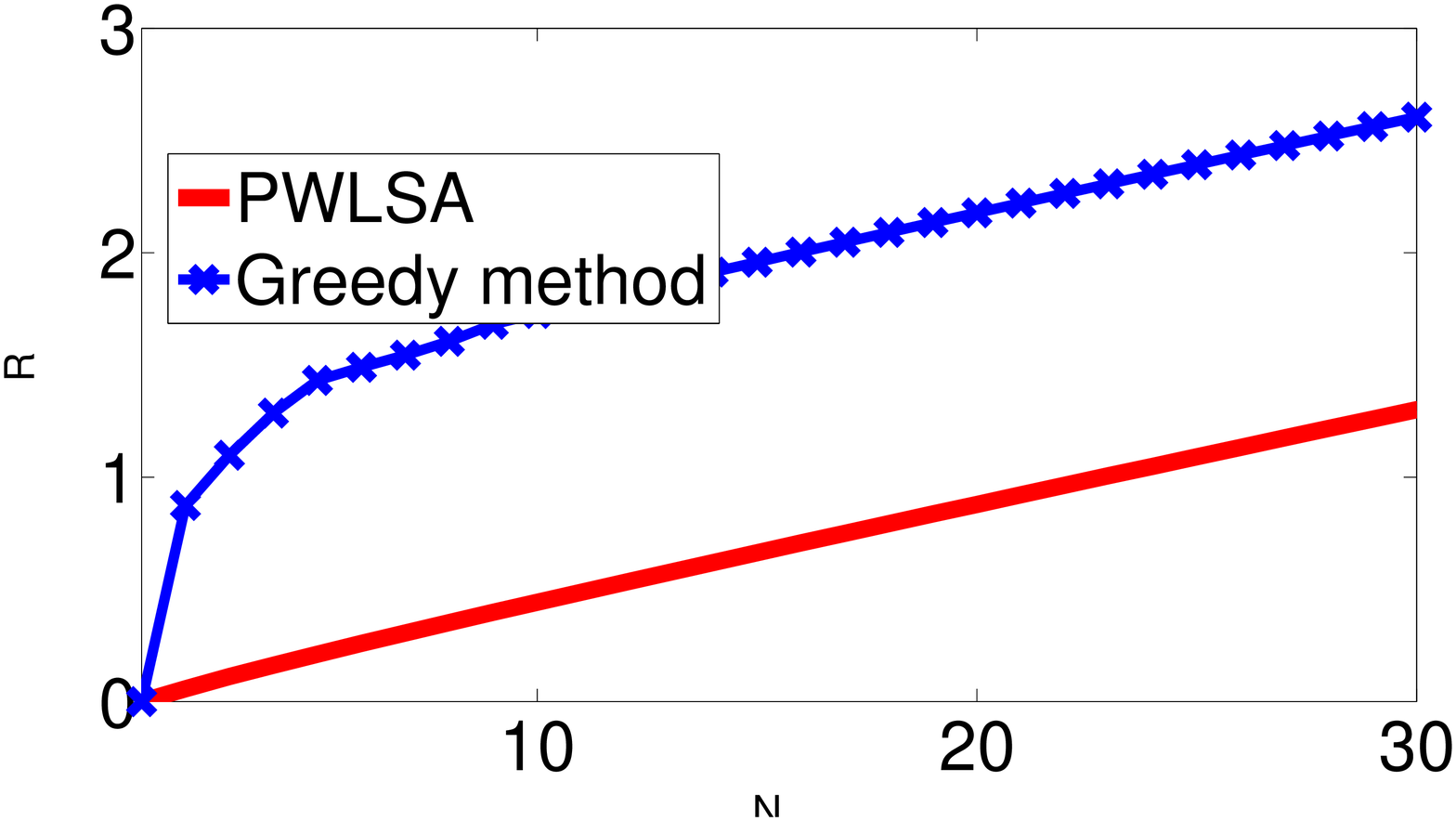}
                \caption{Cumulative regret}
                \label{fig:dyn_avg_regret}
                \end{psfrags}
        \end{subfigure}%
        \begin{subfigure}[b]{0.20\textwidth}
                \begin{psfrags}
                \psfrag{N}[c]{\tiny{Number of Days}}
                \psfrag{D}[c]{\tiny{Price Deviation}}
                \includegraphics[width=1.5in]{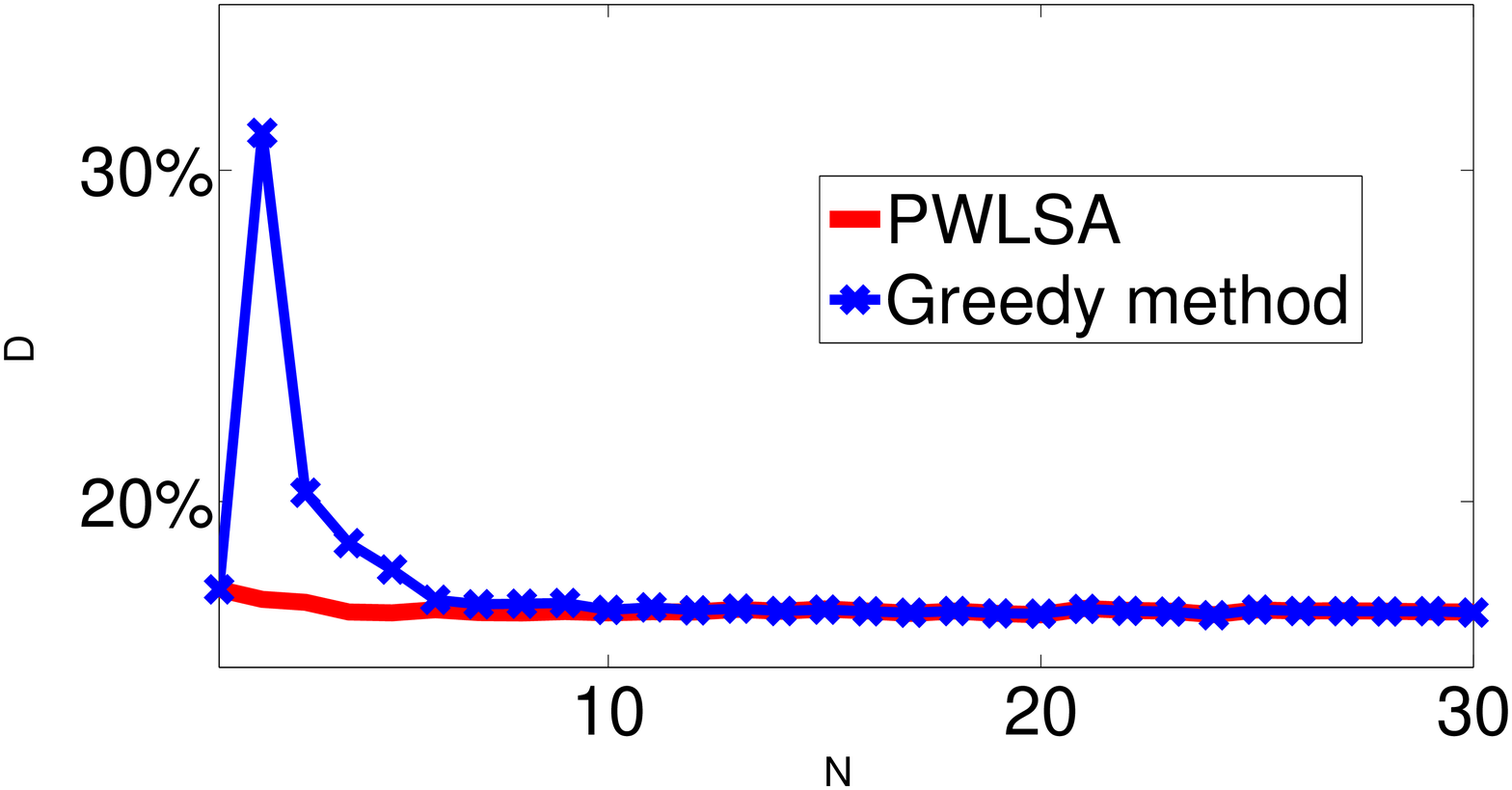}
                \caption{Price convergence}
                \label{fig:dyn_avg_price}
                \end{psfrags}
        \end{subfigure}%
\caption{Average performance comparison of PWLSA and Greedy Method under dynamic unknown demand model, 10,000 Monte Carlo runs}
\label{fig:dyn_avg}
\end{figure}

Fig.~\ref{fig:dyn_avg} shows the average performance comparison of PWLSA and the Greedy Method under dynamic unknown demand model. We can see that PWLSA still outperformed the Greedy Method. According to Fig.~\ref{fig:dyn_avg_regret}, the cumulative regret under PWLSA grew linearly. Intuitively, when a sequence of observation is given, a policy will produce a fixed price or a fixed probability distribution over candidate prices (for randomized policy). However, since the next optimal price is random, there always exists a fixed addition to the expected cumulative regret. Therefore, the linear order achieved by PWLSA is already the best. Fig.~\ref{fig:dyn_avg_price} shows that the error at the very beginning was the cause of the performance gap between PWLSA and the Greedy Method.

\begin{figure}

        \centering
        \begin{subfigure}[b]{0.21\textwidth}
                \begin{psfrags}
                \psfrag{N}[c]{\tiny{Number of Days}}
                \psfrag{R}[c]{\tiny{Price Ratio}}
                \includegraphics[width=1.5in]{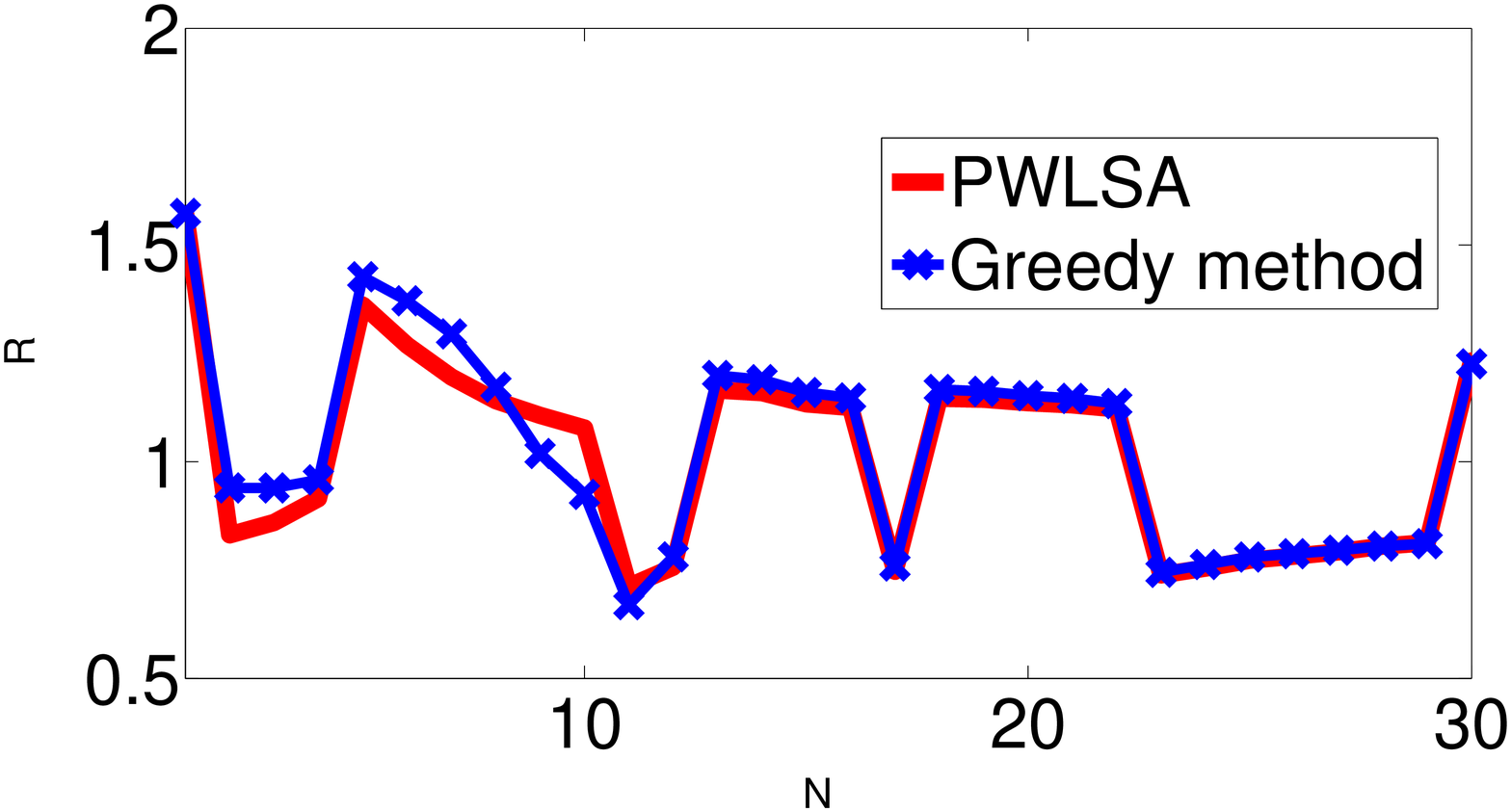}
                \caption{Scenario 1: The two had similar performance}
                \label{fig:dyn_scenario1}
                \end{psfrags}
        \end{subfigure}%
        \begin{subfigure}[b]{0.03\textwidth}
        \end{subfigure}%
        \begin{subfigure}[b]{0.21\textwidth}
                \begin{psfrags}
                \psfrag{N}[c]{\tiny{Number of Days}}
                \psfrag{R}[c]{\tiny{Price Ratio}}
                \includegraphics[width=1.5in]{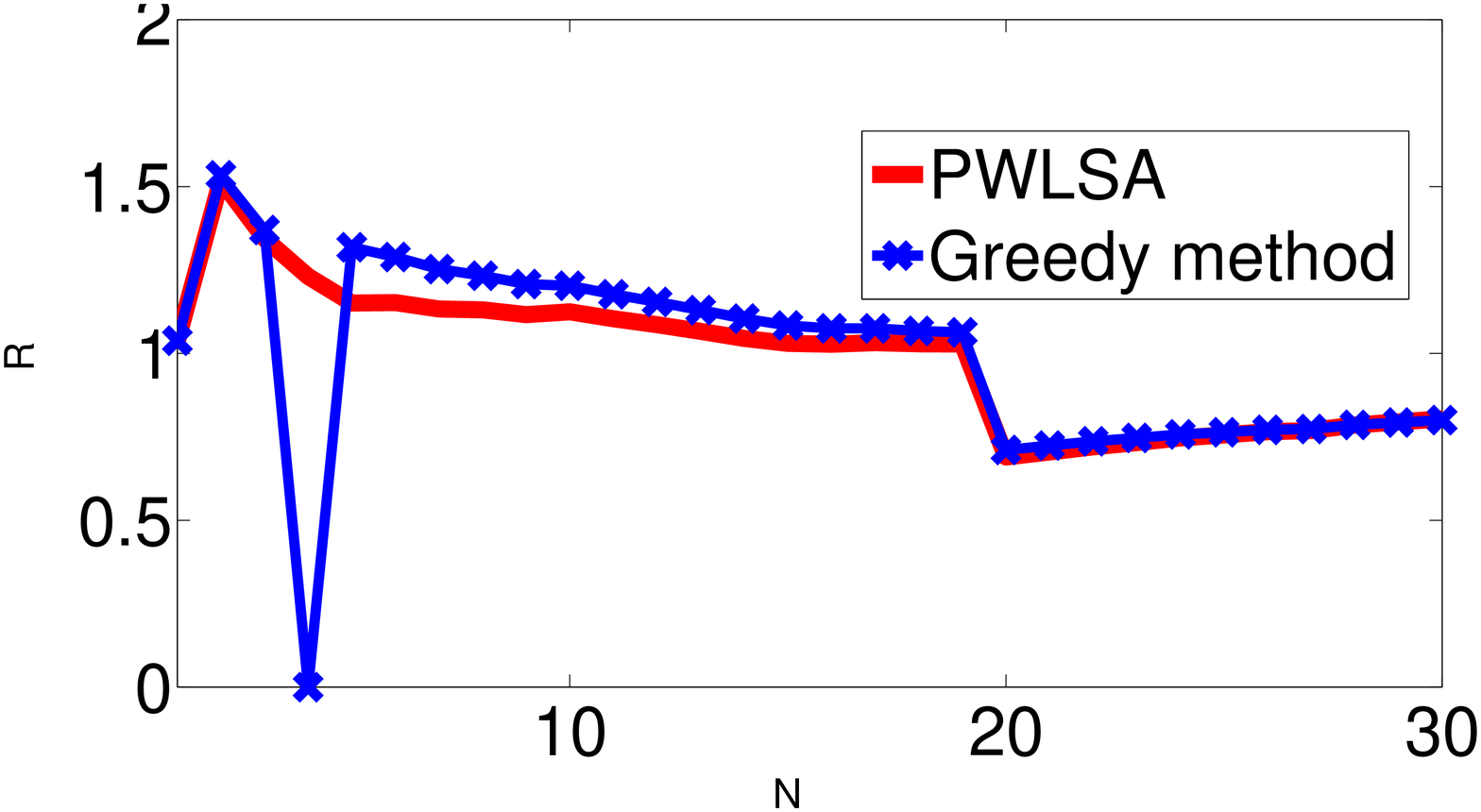}
                \caption{Scenario 2: the Greedy Method failed}
                \label{fig:dyn_scenario2}
                \end{psfrags}
        \end{subfigure}%
\caption{Scenario analysis of PWLSA and Greedy Method under dynamic parameters}
\label{fig:dyn_scenario}
\end{figure}

We also conducted the scenario analysis similar to the static parameter case as shown in Fig.~\ref{fig:dyn_scenario}. Each changing point stands for a incident that the state jumps to the other set. In Fig.~\ref{fig:dyn_scenario1}, we can see that the two policies had similar performance and both tracked the optimal prices well. In few extreme cases, the Greedy Method lost track on the optimal prices wildly at the beginning a few days, as shown in Fig.~\ref{fig:dyn_scenario2}.

\section{Conclusion}

We present in this paper an online learning approach to the dynamic pricing of electricity of a retailer whose customers have price responsive dynamic load with unknown demand function.  We exploit the linear form of the demand function for thermal dynamic load, and cast the problem of online learning as tracking day-ahead dispatch.  This approach leads to a simple learning algorithm with the growth rate of cumulative regret at the order of $\log T$, which is the best rate achievable for any dynamic pricing policies.



\ifCLASSOPTIONcaptionsoff
  \newpage
\fi

\bibliographystyle{IEEEbib}
{
\bibliography{JiaZhaoTong14TSG_final.bbl}
}

\section*{Appendix}

\subsection*{Proof of Lemma~\ref{lemma:square}}
Consider the first order approximation,
\begin{equation}
u(d^{\RT}) - u(d^{\DA}) \approx [\frac{\partial u}{\partial d}(d^{\DA})]^{\mbox{\tiny T}}(d^{\RT} - d^{\DA}).
\end{equation}
By the KKT condition and the definition of real-time price
\begin{equation}
\frac{\partial u}{\partial d}(d^{\DA}) = \frac{\partial c}{\partial p}(p^{\DA}) = \frac{\partial c}{\partial p}(d^{\DA}).
\end{equation}
\begin{equation}
\lambda^{\RT} = \frac{\partial \tilde{c}}{\partial p}(p^{\RT}) = \frac{\partial \tilde{c}}{\partial p}(d^{\RT}) = \frac{\partial c}{\partial p}(d^{\RT}) + \frac{\partial \Delta c}{\partial p}(d^{\RT}),
\end{equation}
where $\Delta c(p) = \tilde{c}(p) - c(p)$. As the generation cost function, $c(p)$ usually takes a quadratic form in practice, $\ie$, $c(p)=\theta p^{\mbox{\tiny T}} p$, where $\theta$ is a scalar. Therefore,
  \begin{equation}
  \label{eq:approxeq}
  \Delta S_{\mbox{\tiny retail}} \approx \theta(d^{\RT} - d^{\DA})^{\mbox{\tiny T}} (d^{\RT} - d^{\DA}) + (\frac{\partial \Delta c}{\partial p}(d^{\RT}))^{\mbox{\tiny T}}(d^{\RT} - d^{\DA}).
  \end{equation}
Usually, the day-ahead cost function $c(p)$ and real-time cost function $\tilde{c}(p)$ have similar shapes and the perturbation $\Delta c$ has small first order derivative. Hence, compared with the first term in (\ref{eq:approxeq}), the second term can be neglected.
\begin{equation}
  \Delta S_{\mbox{\tiny retail}} \approx \theta(d^{\RT} - d^{\DA})^{\mbox{\tiny T}} (d^{\RT} - d^{\DA}).
  \end{equation}

  \hfill $\blacksquare$

\subsection*{Proof of Theorem~\ref{thm:lowerbound}}

First, we reduce the problem to the case that $A$ is known, $d_i^{\DA}$'s are constant and $\Sigma_w$ is a diagonal matrix with the diagonal elements all as $\sigma_w^2$. The minimax rate for this case lower bounds the general case.

For any policy $\mu$, the maximum regret among all possible $b$ is
\begin{equation}
\label{eq:minimax}
L(\mu) = \max_{b} \sum_{t=1}^{T} R_t^{\mu}.
\end{equation}

Assume the parameter $b$ follows a prior distribution,$\gamma_n : \mathscr{N} (\bar{b}, n \sigma^2 \mbox{I})$, where $n$ is a positive integer number and $\mbox{I}$ is an identity matrix. Define the Bayesian cost as (\ref{eq:regret}) and denote the Bayesian estimator of $b$ as $\eta_n$. By the property of joint Gaussian distribution and Sherman-Morrison formula, we can get the minimum Bayesian risk,
\[
R_t^{\eta_n}(\gamma_n) = \mathbb{E}^{\eta_n}||b - \eta_n(d_{1}^{\RT},...,d_{t-1}^{\RT})||_2^2 = \frac{n \sigma^2}{\sigma_w^2 + tn\sigma^2}||\Sigma_w||_2
\]
Then, the cumulative Bayesian risk $\sum_{t=1}^{T} R_t^{\mbox{\tiny Bayes}}(\gamma_n)$ is an increasing function of $n$ and goes to $L(\tilde{\mu}) =  \sum_{t=1}^{T}\frac{1}{t}||\Sigma_w||_2$ as $n$ goes to $\infty$, where $\tilde{\mu}$ is defined in Eq. (\ref{eq:knownA}).

If $\tilde{\mu}$ is not the minimax estimator, there exist some policy $\bar{\mu}$ and $\epsilon >0$, s.t. $L(\bar{\mu}) < L(\tilde{\mu}) - \epsilon$. On the other hand, for $\epsilon >0$, we can find some positive integer $m$, s.t.
\[
L(\tilde{\mu}) - \epsilon <  \sum_{t=1}^{T} R_t^{\eta_m}(\gamma_m).
\]

By the definition (\ref{eq:minimax}), the Bayesian risk of policy $\bar{\mu}$ under the distribution $\gamma_m$ should be less than the maximum cost over all possible values of $b$, $\ie$
\[
\sum_{t=1}^{T} R^{\bar{\mu}}_t (\gamma_m) < L(\bar{\mu}) < L(\tilde{\mu}) - \epsilon <  \sum_{t=1}^{T} R_t^{\eta_m}(\gamma_m),
\]
which contradicts the fact that $\eta_m$ is the Bayesian estimator.


\hfill $\blacksquare$

\subsection*{Proof of Theorem~\ref{thm:logn}}
First, we consider when there is a single day-ahead dispatch level $d^{\DA}$, and $\pi^* = A^{-1}(b-d^{\DA})$. After simplification,
\begin{equation}
\label{eq:sreg}
\begin{array}{l}
\pi_{n+1} - \pi_{n+1}^* \\
= (I-\gamma A)[\Pi_{i=1}^{n-1}(1-\frac{\gamma A}{i+1})](\pi_1 - \pi^*)\\
~~~+ \sum_{k=1}^{n}\{\frac{\gamma}{n} + \sum_{j=k}^{n-1}[\Pi_{i=j+1}^{n-1}(I-\frac{\gamma A}{i+1})] \frac{(I-\gamma A)\gamma}{j(j+1)}\} \omega_k.\\
\end{array}
\end{equation}
For the first term in (\ref{eq:sreg}),
\[
||(I-\gamma A)[\Pi_{i=1}^{n-1}(1-\frac{\gamma A}{i+1})]||_2^2 \le ||(I-\gamma A)||_2^2\Pi_{i=1}^{n-1} ||(I-\frac{\gamma A}{i+1})||_2^2.
\]
Since $(I-\frac{\gamma A}{i+1})^{\mbox{\tiny T}}(1-\frac{\gamma A}{i+1}) = I - \frac{2\gamma A}{i+1} + \frac{\gamma^2 A^2}{(i+1)^2}$, denoting $\lambda_m$ as the minimum eigenvalue of $A$, we have,
\[
||(I-\frac{\gamma A}{i+1})||_2^2 \le I- \frac{2\gamma \lambda_m}{i+1} + \frac{\gamma^2}{(i+1)^2} ||A||_2^2.
\]
Let $C_1 \defeq ||I-\gamma A||_2^2$. Then, since $\gamma \lambda_m > \frac{1}{2}$
\[
\begin{array}{l}
||(I-\gamma A)[\Pi_{i=1}^{n-1}(I-\frac{\gamma A}{i+1})]||_2^2 \\
\le C_1 \Pi_{i=1}^{n-1} (I- \frac{2\gamma \lambda_m}{i+1} + \frac{\gamma^2}{(i+1)^2} ||A||_2^2)= C_2 \frac{1}{n+1}, \\

\end{array}
\]
where $C_2 = C_1 \text{exp}\{\gamma^2||A||_2^2 \}$ doesn't depend on $n$.

For the second term in (\ref{eq:sreg}),
\[
\begin{array}{l}
||[\Pi_{i=j+1}^{n-1}(I-\frac{\gamma A}{i+1})] (I-\gamma A)||_2^2 \\
\le  C_1 \text{exp}\{\sum_{i=j+1}^{n} - \frac{2\gamma \lambda_m}{i+1} + \frac{\gamma^2}{(i+1)^2} ||A||_2^2\} \le C_2 (\frac{j+1}{n+1})^{2\gamma\lambda_m}. \\
\end{array}
\]
Then,
\[
\begin{array}{l}
||\frac{\gamma}{n} + \sum_{j=k}^{n-1}[\Pi_{i=j+1}^{n-1}(I-\frac{\gamma A}{i+1})] \frac{(1-\gamma A)\gamma}{j(j+1)}||_2^2 \\
\le   \{ \frac{\gamma}{n} + \sum_{j=k}^{n-1}||[\Pi_{i=j+1}^{n-1}(I-\frac{\gamma A}{i+1})](I-\gamma A)||_2 \frac{\gamma}{j (j+1)}  \}^2 \\
\le  2 \frac{\gamma^2}{n^2} + 2 \gamma C_2 (\frac{1}{n}) (\frac{1}{n})^{2\gamma\lambda_m - 1} (\frac{1}{k})^{2 - 2\gamma\lambda_m}.
\end{array}
\]
Sum the two terms up,
\[
\begin{array}{l}
\sum_{k=1}^{n-1} ||\frac{\gamma}{n} + \sum_{j=k}^{n-1}[\Pi_{i=j+1}^{n-1}(I-\frac{\gamma A}{i+1})] \frac{(1-\gamma A)\gamma}{j(j+1)}||_2^2 \\
\le 2 \frac{\gamma^2}{n} + 2 \gamma C_2 (\frac{1}{n}) (\frac{1}{n})^{2\gamma\lambda_m - 1} \sum_{k=1}^{n-1} (\frac{1}{k})^{2 - 2\gamma\lambda_m} \le C_3 \frac{1}{n}.\\
\end{array}
\]

Define $M = \max\{||\pi_1 - \pi^*||_2^2, ||\Sigma_\omega||_2^2, ||\Sigma_d||_2^2\}$, we have
\[
\begin{array}{r l}
\sum_{i=1}^{n} L_n & = \mathbb{E} \sum_{i=1}^{n} ||A(\pi_i - \pi^*)||_2^2 \\
& \le \sum_{n=1}^{\mbox{\small T}} ||A||_2^2 [ (C_2 + C_1)\frac{1}{n}] M \le C \log (T).
\end{array}
\]

If $|\mathscr{D}|$ is finite, and we use a separate stochastic approximation to calculate the retail price, then the accumulated regret $\sum_{n=1}^{\mbox{\small T}} R_n \le C|\mathscr{D}| \log (T)$.
\hfill $\blacksquare$

\end{document}

%% file: LTmacros.tex
\newtheorem{theorem}{Theorem}

\newtheorem{lemma}{Lemma}
\newtheorem{definition}{Definition}

\def\beq{\begin{equation}}
\def\eeq{\end{equation}}
\def\bea{\begin{eqnarray}}
\def\eea{\end{eqnarray}}
\def\ba{\begin{array}}
\def\ea{\end{array}}
\def\defeq{{\stackrel{\Delta}{=}}}

\def\eg{{\it e.g., \/}}

\def\ie{{\it i.e.,\ \/}}

\definecolor{bgrd}{rgb}{1,1,1}
\definecolor{grey}{rgb}{0.9,0.9,0.6}
\definecolor{gray}{rgb}{0.5,0.5,0.5}

\makeatletter
\newdimen{\captionwidth}
\long\def\@makecaption#1#2{%
\captionwidth .9\hsize
\vskip 10pt%
\setbox\@tempboxa\hbox{#1: #2}%
  \ifdim \wd\@tempboxa >\captionwidth%
    \setbox\@tempboxa\hbox{#1:\hspace*{.5em}}%
    \hfil\parbox{\captionwidth}{\raggedright\hangindent \wd\@tempboxa%
    \hangafter=1\unhbox\@tempboxa#2}\hfill%
  \else\centerline{\box\@tempboxa}%
  \fi
}
\makeatother

\def\span{\mbox{span}}

\newcommand{\Dmsc}{\mathscr{D}}

\def\Rc{{\cal R}}

%% file: JiaZhaoTong14TSG_final.bbl
\begin{thebibliography}{10}
\providecommand{\url}[1]{#1}
\csname url@samestyle\endcsname
\providecommand{\newblock}{\relax}
\providecommand{\bibinfo}[2]{#2}
\providecommand{\BIBentrySTDinterwordspacing}{\spaceskip=0pt\relax}
\providecommand{\BIBentryALTinterwordstretchfactor}{4}
\providecommand{\BIBentryALTinterwordspacing}{\spaceskip=\fontdimen2\font plus
\BIBentryALTinterwordstretchfactor\fontdimen3\font minus
  \fontdimen4\font\relax}
\providecommand{\BIBforeignlanguage}[2]{{%
\expandafter\ifx\csname l@#1\endcsname\relax
\typeout{** WARNING: IEEEtran.bst: No hyphenation pattern has been}%
\typeout{** loaded for the language `#1'. Using the pattern for}%
\typeout{** the default language instead.}%
\else
\language=\csname l@#1\endcsname
\fi
#2}}
\providecommand{\BIBdecl}{\relax}
\BIBdecl

\bibitem{Borenstein&etc:02}
S.~Borenstein, M.~Jaske, and A.~Rosenfeld, ``{Dynamic Pricing, Advanced
  Metering, and Demand Response in Electricity Markets},'' \emph{Recent Work,
  Center for the Study of Energy Markets, University of California Energy
  Institute, UC Berkeley}, 2002.

\bibitem{Hopper&Goldman&Neenan:05EJ}
N.~Hopper, C.~Goldman, and B.~Neenan, ``{Demand Response from Day-Ahead Hourly
  Pricing for Large Customers},'' \emph{Electricity Journal}, no.~02, pp.
  52--63, 2006.

\bibitem{Borenstein:05}
S.~Borenstein, ``{The long run efficiency of real-time electricity pricing},''
  \emph{The Energy Journal}, vol.~26, no.~3, 2005.

\bibitem{Carrion&Etal:07TPS}
M.~Carrion, A.~Conejo, and J.~Arroyo, ``{Forward Contracting and Selling Price
  Determination for a Retailer},'' \emph{IEEE Transactions on Power Systems},
  vol.~32, no.~4, November 2007.

\bibitem{Conejo&Etal:08TPS}
A.~Conejo, R.~Garcia-Bertrand, M.~Carrion, A.~Caballero, and A.~Andres,
  ``{Optimal Involvement in Futures Markets of a Power Producer},'' \emph{IEEE
  Transactions on Power Systems}, vol.~23, no.~2, May 2008.

\bibitem{Yang&Etal:12TPS}
P.~Yang, G.~Tang, and A.~Nehorai, ``{A Game-Theoretic Approach for Optimal
  Time-of-Usa Electricity Pricing},'' \emph{IEEE Transactions on Power
  Systems}, 2012.

\bibitem{Jia&Tong:12Allerton}
L.~Jia and L.~Tong, ``Optimal pricing for residential demand response: A
  stochastic optimization approach,'' in \emph{2012 Allerton Conference on
  Communication, Control and Computing}, Oct. 2012.

\bibitem{Jia&Tong:13CDC}
------, ``Day ahead dynamic pricing for demand response in dynamic
  environments,'' in \emph{52nd IEEE Conference on Decision and Control}, Dec.
  2013.

\bibitem{Lai&Robbins:85AAM}
T.~L. Lai and H.~Robbins, ``{Asymptotically efficient adaptive allocation
  rules},'' \emph{Advanced and Applied Mathematics}, vol.~6, no.~1, pp. 4--22,
  1985.

\bibitem{Agrawal:95SIAM}
R.~Agrawal, ``{The Continuum-Armed Bandit Problem},'' \emph{SIAM J. Control and
  Optimization}, vol.~33, no.~6, pp. 1926--1951, 1995.

\bibitem{Kleinberg:04}
R.~Kleinberg, ``{Nearly Tight Bounds for the Continuum-Armed Bandit Problem},''
  \emph{Advances in Neural Information Processing Systems}, pp. 697--740, 2004.

\bibitem{Auer&etal:07}
P.~Auer, R.~Ortner, and C.~Szepesvari, ``{Improved Rates for the Stochastic
  Continuum-Armed Bandit Problem},'' \emph{Lecture Notes in Computer Science},
  vol. 4539, pp. 454--468, 2007.

\bibitem{Cope:09}
E.~W. Cope, ``{Regret and Convergence Bounds for a Class of Continuum-Armed
  Bandit Problems},'' \emph{IEEE Transactions on Automatic Control}, vol.~54,
  no.~6, Jan. 2009.

\bibitem{KleinbergLeighton:03}
R.~Kleinberg and T.~Leighton, ``{The value of knowing a demand curve: bounds on
  regret for online posted-price auctions},'' in \emph{Proc. 44th IEEE
  Symposium on Foundations of Computer Science (FOCS)}, 2003.

\bibitem{RusmevichientongTsitsiklis:10}
P.~Rusmevichientong and J.~N. Tsitsiklis, ``{Linearly Parameterized Bandits},''
  \emph{Mathematics of Operations Research}, vol.~35, no.~2, pp. 395--411,
  2010.

\bibitem{BroderRusmevichientong:12}
J.~Broder and P.~Rusmevichientong, ``{Dynamic Pricing under a General
  Parametric Choice Model},'' \emph{Operations Research}, vol.~60, no.~4, pp.
  965--980, 2012.

\bibitem{Harrison&etal:SC11}
J.~M. Harrison, N.~B. Keskin, and A.~Zeevi, ``{Bayesian Dynamic Pricing
  Policies: Learning and Earning Under a Binary Prior Distribution},''
  \emph{Management Science}, pp. 1--17, Oct. 2011.

\bibitem{Zhai&etal:Asilomar11}
Y.~Zhai, P.~Tehrani, L.~Li, J.~Zhao, and Q.~Zhao, ``{Dynamic Pricing under
  Binary Demand Uncertainty: A Multi-Armed Bandit with Correlated Arms},'' in
  \emph{Proc. of the 45th IEEE Asilomar Conference on Signals, Systems, and
  Computers}, Nov. 2011.

\bibitem{Lai&Robbins:79AS}
T.~L. Lai and H.~Robbins, ``{Adaptive design and stochastic approximation},''
  \emph{The Annals of Statistics}, vol.~7, no.~6, pp. 1196--1221, 1979.

\bibitem{Lai&Robbins:82AAM}
------, ``{Iterated Least Squares in Multiperiod Control},'' \emph{Advanced and
  Applied Mathematics}, vol.~3, pp. 50--73, 1982.

\bibitem{Bersimas&Perakis:06}
D.~Bertsimas and G.~Perakis, ``{Dynamic Pricing: A Learning approach},''
  \emph{Mathematical and Computational Models for Congestion Charging}, pp.
  45--80, 2006.

\bibitem{Lobo&Boyd:03Informs}
M.~Lobo and S.~Boyd, ``{Pricing and learning with uncertain demand},'' in
  \emph{INFORMS Revenue Management Conference}, Columbia University, 2003.

\bibitem{Garcia&etal:05OP}
E.~C.-N. Alfredo~Garcia and J.~Reitzes, ``{Dynamic Pricing and Learning in
  Electricity Markets},'' \emph{Operation Research}, vol.~53, no.~2, pp.
  231--241, 2005.

\bibitem{RahimiKianSadeghiThomas05PES}
A.~Rahimi-Kian, B.~Sadeghi, and R.~J. Thomas, ``{Q-learning based
  supplier-agents for electricity markets},'' in \emph{IEEE Power Engineering
  Society General Meeting}, 2005.

\bibitem{QiuPeetersDeconinck09ISAP}
Z.~Qiu, E.~Peeters, and G.~Deconinck, ``{Comparison of Two Learning Algorithms
  in Modelling the Generator's Learning Abilities},'' in \emph{15th
  International Conference on Intelligent System Applications to Power
  Systems}, 2009.

\bibitem{Pinto&Etal:11ISAP}
T.~Pinto, Z.~Vale, F.~Rodrigues, and I.~Praca, ``{Cost dependent strategy for
  electricity markets bidding based on adaptive reinforcement learning},'' in
  \emph{the 16th International Conference onIntelligent System Application to
  Power Systems}, 2011.

\bibitem{Taylor&Mathieu:13STPS}
J.~A. Taylor and J.~L. Mathieu, ``{Index Policies for Demand Response},''
  \emph{IEEE Transactions on Power System}, vol.~PP, no.~99, pp. 1--9, 2013.

\bibitem{O'Neill&etal:10SGC}
A.~G. Daniel~O'Neill, Marco~Levorato and U.~Mitra, ``{Residential Demand
  Response Using Reinforcement Learning},'' in \emph{Proceedings of 2010 First
  IEEE International Conference on Smart Grid Communications}, Oct. 2010, pp.
  409--414.

\bibitem{EnergySurvey}
Residential Energy Consumption Survey (RECS),
  http://www.eia.gov/consumption/residential/.

\bibitem{Bargiotas&Birddwell:88ITPD}
D.~Bargiotas and J.~Birddwell, ``{Residential air conditioner dynamic model for
  direct load control},'' \emph{IEEE Transactions on Power Delivery}, vol.~3,
  no.~4, pp. 2119--2126, Oct. 1988.

\bibitem{Yu&etal:13TSG}
Z.~Yu, L.~McLaughlin, L.~Jia, M.~C. Murphy-Hoye, A.~Pratt, and L.~Tong,
  ``{Modeling and Stochastic Control for Home Energy Management},'' \emph{IEEE
  Transactions on Smart Grid}, 2013.

\bibitem{Anderson&Taylor:76}
T.~Anderson and J.~Taylor, ``{Some experimental results on the statistical
  properties of least squares estimates in control prob- lems},''
  \emph{Econometrica}, vol.~44, p. 1289¨C1302, 1976.

\end{thebibliography}
